\documentclass[a4paper,12pt]{article}
\usepackage{comment}
\usepackage{cite}
\usepackage{amsmath}
\usepackage{amssymb}
\usepackage{amsfonts}
\usepackage[T1]{fontenc}
\usepackage[utf8]{inputenc}
\usepackage{graphicx}
\usepackage{fancyhdr}
\usepackage{float}
\usepackage{xcolor}
\usepackage{authblk}
\usepackage{mathrsfs}
\usepackage{empheq}
\usepackage[hyphens]{url}
\usepackage{hyperref} 
\usepackage[]{breakurl}
\usepackage[]{amsthm}
\usepackage{tikz}

\newtheorem{theorem}{Theorem}

\usepackage{geometry}
\begin{document}

\title{An integral representation of Catalan numbers using Malmst\'en's formula}

\author[$\dagger$]{Jean-Christophe {\sc Pain}$^{1,2,}$\footnote{jean-christophe.pain@cea.fr}\\
\small
$^1$CEA, DAM, DIF, F-91297 Arpajon, France\\
$^2$Universit\'e Paris-Saclay, CEA, Laboratoire Mati\`ere en Conditions Extr\^emes,\\ 
91680 Bruy\`eres-le-Ch\^atel, France
}

\maketitle

\begin{abstract}
In this article, we propose an integral expression of the Catalan numbers, based on Malmst\'en's definite-integral representation of $\ln\left[\Gamma(x)\right]$, $\Gamma$ being the usual Gamma function. The obtained expression is likely to yield new summations involving Catalan numbers or central binomial coefficients.
\end{abstract}

\section{Introduction}

The Catalan numbers, defined as
\begin{equation*}
C_n=\frac{1}{n+1}\binom{2n}{n}=~_2F_1\left[\begin{array}{cc}
1-n,-n\\
2
\end{array};1\right],
\end{equation*}
where $_2F_1$ represents the Gauss hypergeometric function, are encountered in many problems of combinatorics \cite{Stanley2015,Koshy2008}. $C_n$ represents, in particular, the number of expressions containing $n$ pairs of correctly-matched parentheses. For instance, on the $n=3$ case, the different possibilities are ((())), ()(()), ()()(), (())() and (()()), and therefore $C_3=5$. $C_n$ represents also the number of ways a convex polygon of $n+2$ sides can split into triangles by connecting vertices:
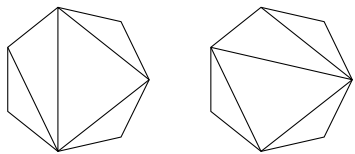
\def\hepta{\draw(A) -- (B) -- (C) -- (D) -- (E) -- (F) -- (G) -- cycle;}


\newcommand{\slice}[1]{%
    \hepta
    \draw \foreach \x/\y in {#1} {(\x)--(\y)};

}
\begin{center}
\begin{tikzpicture}
    \coordinate (A) at (-0.76,1.54);
    \coordinate (B) at (-0.76,0.69);
    \coordinate (C) at (-0.10,0.16);
    \coordinate (D) at (0.73,0.35);
    \coordinate (E) at (1.1,1.11);
    \coordinate (F) at (0.73,1.88);
    \coordinate (G) at (-0.10,2.07);

\matrix[column sep=0.8cm,row sep=0.5cm]
{
    \slice{A/C,C/E,E/G,C/G}&
    \slice{A/C,C/E,E/G,A/E}&
\\
};
\end{tikzpicture}
\end{center}

For a heptagon, the number of such diagrams is $C_5=42$. Although several integral representations of Catalan numbers have been obtained (see for instance the non-exhaustive list of references \cite{Penson2001,Dana2012,Shi2015,Qi2017,Guo2023}), this is still an active field of research \cite{Qi2024}. For instance, one can write \cite{Penson2001}:
\begin{equation*}
C_{n}=\frac{2}{\pi}4^{n}\int_{-1}^{1}t^{2n}\sqrt{1-t^{2}}\,dt,
\end{equation*}
or equivalently
\begin{equation*}
C_{n}=\frac{4^{n+2}}{\pi}\int_{0}^{\infty}\frac{\sqrt{t}}{(4t+1)^{n+2}}\,dt,
\end{equation*}
As another example, in 2015, Xiao-Ting Shi, Fang-Fang Liu and Feng Qi obtained a new expression using the Binet formula integral \cite{Temme1996}. The Binet \cite{Whittaker1990,Sasvari1999} formula reads: 
\begin{equation}\label{bin}
    \Gamma(x+1)=\left(\frac{x}{e}\right)^x~\sqrt{2\pi x}~e^{\theta(x)},
\end{equation}
where
\begin{equation*}
\Gamma(x+1)=\int_0^{\infty}\,t^x\,e^{-t}\,\mathrm{dt}
\end{equation*}
is the usual Gamma function and
\begin{equation*}
    \theta(x)=\int_0^{\infty}\left(\frac{1}{e^t-1}-\frac{1}{t}+\frac{1}{2}\right)~\frac{e^{-xt}}{t}~\mathrm{d}t.
\end{equation*}
Equation (\ref{bin}) can be recast into
\begin{equation*}
    \ln\left[\Gamma(x+1)\right]=x\ln x-x+\frac{1}{2}\ln(2\pi x)+\theta(x),
\end{equation*}
and Xiao-Ting Shi, Fang-Fang Liu and Feng Qi derived the following expression:
\begin{equation*}\label{feng}
C_n=\frac{e^{3/2}\,4^n(n+1/2)^n}{\sqrt{\pi}\,(n+2)^{n+3/2}}\,\exp\left[\int_0^{\infty}\left(\frac{1}{e^t-1}-\frac{1}{t}+\frac{1}{2}\right)\frac{\left(e^{-t/2}-e^{-2t}\right)}{t}\,e^{-xt}\,\mathrm{d}t\right].    
\end{equation*}

In a different framework, we recently obtained two integral representations \cite{Pain2024b} using the Touchard identity \cite{Touchard1928}, and a variant of the latter published by Callan \cite{Callan2013}. We also proposed integral representations of the Glaisher-Kinkelin constant $A$ \cite{Glaisher1878,Pain2024a} using integral representations of the logarithm of the Gamma function \cite{Boros2004} together with the following integral \cite{Glaisher1878,Choi1997}:
\begin{equation}\label{gla}
    \int_0^{1/2}\ln\left[\Gamma(x+1)\right]~\mathrm{d}x=-\frac{1}{2}-\frac{7}{24}\ln 2+\frac{1}{4}\ln\pi+\frac{3}{2}\ln A.
\end{equation}

Integral representations of Catalan numbers can be of great interest to deduce sum rules for the latter. We have in particular \cite{Stewart2022}:
\begin{equation*}
\sum_{n=0}^{\infty}\frac{1}{(2n+1)\,64^n}\,C_{2n}\,C_n=\frac{8\sqrt{2}}{3\pi},   
\end{equation*}
or also
\begin{equation*}
\sum_{n=0}^{\infty}\frac{C_{2n}\,C_n}{64^n}=\frac{4}{\pi}\,\ln(3+2\sqrt{2})-\frac{8\sqrt{2}}{3\pi}.
\end{equation*}
In section \ref{sec2}, we obtain, using the Malmst\'en integral for the logarithm of the Gamma function, an integral representation for the Catalan numbers. The resulting formula is different from the one in Eq. (\ref{feng}).  

\section{Integral representation deduced from the Malmst\'en formula}\label{sec2} 

\begin{theorem}

Let $n$ be a positive integer and $C_n$ the $n^{th}$ Catalan number. Then the integral representation
\begin{equation*}
    C_n=\frac{4^n}{\sqrt{\pi}}\,\exp\left\{\int_0^{\infty}\left[\frac{(e^{3t/2}-1)}{(e^{t}-1)}e^{-nt}-\frac{3}{2}\right]\frac{e^{-t}}{t}dt\right\}
\end{equation*}
holds true.

\end{theorem}

\begin{proof}

Using the alternative representation of Catalan numbers:

\begin{equation*}
    C_n=\frac{4^n\Gamma(n+1/2)}{\sqrt{\pi}\,\Gamma(n+2)},
\end{equation*}
one gets, for $n\geq 1$:
\begin{equation}\label{lncn}
    \ln\left(C_n\right)=(2\ln 2)\,n-\frac{\ln\pi}{2} +\ln\left[\Gamma\left(n+\frac{1}{2}\right)\right]-\ln\left[\Gamma(n+2)\right]
\end{equation}
The Malmst\'en formula reads \cite{Erdelyi1981,Pain2024b}:
\begin{equation*}
    \ln\left[\Gamma(x+1)\right]=\int_0^{\infty}\left[x-\frac{(1-e^{-xt})}{(1-e^{-t})}\right]\frac{e^{-t}}{t}dt,
\end{equation*}
yielding in particular, and for our purpose:
\begin{equation}\label{a1}
    \ln\left[\Gamma\left(x+\frac{1}{2}\right)\right]=\int_0^{\infty}\left[x-\frac{1}{2}-\frac{(1-e^{-\left(x-\frac{1}{2}\right)t})}{(1-e^{-t})}\right]\frac{e^{-t}}{t}dt,
\end{equation}
as well as
\begin{equation}\label{a2}
    \ln\left[\Gamma(x+2)\right]=\int_0^{\infty}\left[x+1-\frac{(1-e^{-(x+1)t})}{(1-e^{-t})}\right]\frac{e^{-t}}{t}dt.
\end{equation}
Using Eqs. (\ref{a1}) and (\ref{a2}), one gets

\begin{equation*}
     \ln\left[\Gamma\left(x+\frac{1}{2}\right)\right]-\ln\left[\Gamma(x+2)\right]=\int_0^{\infty}\left[\frac{(e^{-t}-e^{t/2})}{(e^{-t}-1)}e^{-tx}-\frac{3}{2}\right]\frac{e^{-t}}{t}dt.   
\end{equation*}
and Eq. (\ref{lncn}) becomes
\begin{equation*}
    \ln\left(C_n\right)=(2\ln 2)\,n-\frac{\ln\pi}{2}+\int_0^{\infty}\left[\frac{(e^{-t}-e^{t/2})}{(e^{-t}-1)}e^{-nt}-\frac{3}{2}\right]\frac{e^{-t}}{t}dt.   
\end{equation*}
and thus
\begin{equation*}
    C_n=\frac{4^n}{\sqrt{\pi}}\,\exp\left\{\int_0^{\infty}\left[\frac{(e^{-t}-e^{t/2})}{(e^{t}-1)}e^{-nt}-\frac{3}{2}\right]\frac{e^{-t}}{t}dt\right\}.   
\end{equation*}
or also
\begin{equation*}
    C_n=\frac{4^n}{\sqrt{\pi}}\,\exp\left\{\int_0^{\infty}\left[\frac{(e^{3t/2}-1)}{(e^{t}-1)}e^{-nt}-\frac{3}{2}\right]\frac{e^{-t}}{t}dt\right\}.   
\end{equation*}
which completes the proof.

\end{proof}

\section{Conclusion}\label{sec3}

In this work, an integral representation of the Catalan numbers was obtained, thanks to an integral of $\ln\left[\Gamma(x+1)\right]$ due to Malmst\'en. It is hoped that the new formula, {\it a priori} different from the existing ones, in the sense that it can not be easily deduced from them, will stimulate the derivation of identities involving (products of) central binomial coefficients and Catalan numbers. We also plan to apply the techniques presented here to generalized Catalan numbers, entering certain matrix models related to hypergraphs \cite{Kahkeshani2013,Gunnells2021}.


\begin{thebibliography}{99}

\bibitem{Stanley2015} R. P. Stanley, {\it Catalan numbers}, Cambridge University Press, 2015.

\bibitem{Koshy2008} T. Koshy, {\it Catalan numbers with applications}, Oxford Academic, New York, 2008. 

\bibitem{Penson2001} K. A. Penson and J.-M. Sixdeniers, Integral representations of Catalan and related numbers, {\it J. Integer Seq.} {\bf 4}, Art. 01.2.5 (2001).

\bibitem{Dana2012} T. Dana-Picard and D. G. Zeitoun, Parametric improper integrals, Wallis formula and Catalan numbers, {\it Internat. J. Math. Ed. Sci. Tech.} {\bf 43}, 515-520 (2012).

\bibitem{Shi2015} Xiao-Ting Shi, Fang-Fang Liu, Feng Qi, An integral representation of the Catalan
numbers, {\it Glob. J. Math. Anal.} {\bf 3}, 130-133 (2015).

\bibitem{Qi2017} Fen Qi and Bai-Ni Guo, Integral representations of the Catalan numbers and their applications, {\it Mathematics} {\bf 5}, 40 (2017). 

\bibitem{Guo2023} Bai-Ni Guo and Dongkyu Lim, Integral representations of Catalan numbers and sums involving central binomial coefficients, {\it Integers} {\bf 23}, \#A34 (2023).

\bibitem{Qi2024} Feng Qi, Some properties of the Catalan numbers, {\it Anal. Discrete Math.} {\bf 19}, 10 pp. (2024).

\bibitem{Temme1996} N. M. Temme, {\it Special functions: an introduction to classical functions of mathematical physics}, John Wiley \& Sons, Inc., New York, 1996.

\bibitem{Whittaker1990} E. T. Whittaker and G. N. Watson, {\it A course in modern analysis}, 4th ed. Cambridge, England: Cambridge University Press, 1990.

\bibitem{Sasvari1999} Z. Sasvari, An elementary proof of Binet's formula for the Gamma function, {\it Amer. Math. Mon.} {\bf 106}, 156-158 (1999).

\bibitem{Pain2024b} J.-C. Pain, Integral representations of Catalan numbers using Touchard-like identities, \url{https://www.arxiv.org/abs/2408.01201} (2024).

\bibitem{Touchard1928} J. Touchard, Sur certaines \'equations fonctionnelles, in: {\it Proc. Int. Math. Congress, Toronto (1924)} {\bf 1}, 465-472 (1928). In French.

\bibitem{Callan2013} D. Callan, A variant of Touchard's Catalan number identity, \url{https://arxiv.org/abs/1204.5704} (2024).

\bibitem{Glaisher1878} J. W. L. Glaisher, On the Product $1^1.2^2.3^3...n^n$, {\it Messenger Math.} {\bf 7}, 43-47, 1878.

\bibitem{Pain2024a} J.-C. Pain, Two integral representations for the logarithm of the Glaisher-Kinkelin constant, \url{https://arxiv.org/abs/2405.05264} (2024).

\bibitem{Boros2004} G. Boros and V. Moll, The expansion of the loggamma function, §10.6 in {\it Irresistible integrals: symbolics, analysis and experiments in the evaluation of integrals}, Cambridge, England: Cambridge University Press, pp. 201-203, 2004.

\bibitem{Choi1997} J. Choi and C. Nash, Integral representations of the Kinkelin's constant $A$, {\it Math. Japon.} {\bf 45}, 223–230 (1997).

\bibitem{Stewart2022} S. M. Stewart, A simple proof and some applications of an integral representation for the Catalan numbers, {\it Applied Mathematics E-Notes} {\bf 22}, 637-643 (2022).

\bibitem{Erdelyi1981} A. Erd\'elyi, W. Magnus, F. Oberhettinger and F. G. Tricomi, {\it Higher transcendental functions}, Vol. 1. New York: Krieger, pp. 20-21, 1981.

\bibitem{Gunnells2021} P. E. Gunnells, Generalized Catalan numbers from hypergraphs, \url{https://arxiv.org/abs/2102.05121}

\bibitem{Kahkeshani2013} R. Kahkeshani, A generalization of the Catalan numbers, {\it J. Integer Seq.} {\bf 16}, Article 13.6.8 (2013).

\end{thebibliography}
\end{document}